\documentclass{article}

\usepackage{amsmath}
\usepackage{amssymb}

\newtheorem{thm}{Theorem}

\newtheorem{lem}{Lemma}

\newtheorem{ques}{Question}

\title{Perfect squares have at most five divisors close to its square root}
\author{Tsz Ho Chan}
\date{}

\begin{document}
\maketitle

\begin{abstract}
In this paper, we consider a conjecture of Erd\H{o}s and Rosenfeld when the number is a perfect square. In particular, we show that every perfect square $n$ can have at most five divisors between $\sqrt{n} - c \sqrt[4]{n}$ and $\sqrt{n} + c \sqrt[4]{n}$.
\end{abstract}

\section{Introduction and main result}
In [\ref{ER}], Erd\H{o}s and Rosenfeld considered the differences between the divisors of a positive integer $n$. They exhibited infinitely many integers with four ``small" differences and posed the question that any positive integer can have at most a bounded number of ``small" differences. Specifically, they asked
\begin{ques} \label{ques1}
Is there an absolute constant $K$, so that for every $c$, the number of divisors of $n$ between $\sqrt{n}$ and $\sqrt{n} + c \sqrt[4]{n}$ is at most $K$ for $n > n_0(c)$?
\end{ques}
In this paper, we answer the above question when $n$ is a perfect square. In particular, we have
\begin{thm} \label{thm1}
For every $c \ge 3$, any perfect square $n$ can have at most five divisors between $\sqrt{n} - c \sqrt[4]{n}$ and $\sqrt{n} + c \sqrt[4]{n}$ for $n > e^{C c^6 (\log c)^5}$ where $C$ is some sufficiently large constant independent of $c$.
\end{thm}
This answers question \ref{ques1} for perfect squares with $K = 3$. Based on the proof of Theorem \ref{thm1}, every example where a perfect square $n$ has three divisors between $\sqrt{n}$ and $\sqrt{n} + c \sqrt[4]{n}$ comes from solutions to Pell equations. For example, consider the Pell equation $X^2 - 2 Y^2 = 2$. It has solutions $(X_k, Y_k)$ generated by $X_k + \sqrt{2} Y_k = (3 + 2\sqrt{2})^k (2 + \sqrt{2})$. Then one can verify that $X_k$ are even, $Y_k$ are odd and $(X_k - 2)(X_k + 2) = 2 (Y_k - 1) (Y_k + 1)$. Now consider the integers $n = (X_k - 2)^2 (X_k + 2)^2 = 4 (Y_k - 1)^2 (Y_k + 1)^2$. It has divisors $(X_k - 2) (X_k + 2)$, $(X_k + 2)^2$, $2 (Y_k + 1)^2$ that are between $\sqrt{n}$ and $\sqrt{n} + 5 \sqrt[4]{n}$. This shows that $K = 3$ is the best possible constant for question \ref{ques1} to be true with perfect squares.

\section{Initial transformation}
Suppose $N^2 = (N - d_1)(N + e_1) = (N - d_2)(N + e_2) = ... = (N - d_r)(N + e_r)$ where $N, N - d_i, N + e_j$ are all the divisors of $N^2$ that lie in $[N - c N^{1/2}, N + c N^{1/2}]$ where $1 \le d_1 < d_2 < ... < d_r \le c N^{1/2}$ and $1 \le e_1 < e_2 < ... < e_r \le c N^{1/2}$ are positive integers. Observe that $N^2 = (N - d_i)(N + e_i)$ which gives $e_i d_i = (e_i - d_i) N$. So we must have $e_i > d_i$ and say $e_i = d_i + l_i$ for some positive integer $l_i$. Hence $(d_i + l_i) d_i = l_i N$ which gives $d_i^2 + l_i d_i = l_i N$. Multiply both sides by four and add $l_i^2 + 4N^2$ to both sides, we have $(2d_i + l_i)^2 + (2N)^2 = (2N + l_i)^2$. Also from $(d_i + l_i) d_i = l_i N$, we have
\begin{equation} \label{restrict}
l_i = \frac{d_i^2}{N - d_i} \le \frac{ c^2 N }{N - c N^{1/2}} \le 2 c^2
\end{equation}
for $N \ge 4c^2$.

\section{Pythagorean triples}
Thus we have a Pythagorean triple $2d_i + l_i$, $2N$, $2N+l_i$. It is well-known that all the solutions to the Pythagorean equation are parametrized by $\lambda (u^2 - v^2)$, $\lambda (2uv)$, $\lambda (u^2 + v^2)$ for some positive integers $\lambda$ and $u > v$.

\bigskip

Case 1: $2d_i + l_i = \lambda_i (u_i^2 - v_i^2)$, $2N = \lambda_i (2u_i v_i)$ and $2N + l_i = \lambda_i (u_i^2 + v_i^2)$.
Subtracting the last two equations, we have $l_i = \lambda_i (u_i - v_i)^2 \le 2 c^2$ by (\ref{restrict}). Hence $u_i - v_i \le \sqrt{2} c$, $\lambda_i \le 2 c^2$ and $2N$ can be written as $2 \lambda_i u_i v_i$. By adding or subtracting the three equations, we also have $2(N - d_i) = 2\lambda_i v_i^2$ and $2(N + e_i) = 2(N + d_i + l_i) = 2\lambda_i u_i^2$.

\bigskip

Case 2: $2N = \lambda_i (u_i^2 - v_i^2)$, $2d_i + l_i = \lambda_i (2u_i v_i)$ and $2N + l_i = \lambda_i (u_i^2 + v_i^2)$. Subtracting the first and the last equations, we have $l_i = 2 \lambda_i v_i^2 \le 2 c^2$ by (\ref{restrict}). Hence $v_i \le c$, $\lambda_i \le c^2$ and $2N$ can be written as $\lambda_i (u_i - v_i) (u_i + v_i)$. By adding or subtracting the three equations, we also have $2(N - d_i) = \lambda_i (u_i - v_i)^2$ and $2(N + e_i) = 2(N + d_i + l_i) = \lambda_i (u_i + v_i)^2$.

\bigskip

In either case, $2N = \mu_i x_i y_i$, $2(N-d_i) = \mu_i x_i^2$ and $2(N+e_i) = \mu_i y_i^2$ with $1 \le y_i - x_i \le 2c$, $\mu_i = \lambda_i$ or $2 \lambda_i$ and $\mu_i \le 4 c^2$.

\section{Almost squares} \label{almostsquare}

Now we claim that the $\mu_i$ are distinct if $N > 32 c^6$. Suppose not, say $\mu_i = \mu_j$ for some $1 \le i < j \le r$. Then $\mu_i x_i y_i = 2N = \mu_j x_j y_j$ implies $x_i y_i = x_j y_j = \frac{2N}{\mu_i}$. Numbers like $\frac{2N}{\mu_i}$ that can be factored as $x_i y_i$ and $x_j y_j$ with $x_i$ close to $y_i$ and $x_j$ close to $y_j$ are called almost squares of type 2 and have been studied by the author in [\ref{C1}], [\ref{C2}] and [\ref{C3}] for example. If $x_i = x_j$, then $2(N - d_i) = \mu_i x_i^2 = \mu_j x_j^2 = 2(N - d_j)$ which contradicts $d_i < d_j$. Without loss of generality, assume $x_i < x_j$. Then we must have $x_i < x_j < y_j < y_i$. Let $y_j = n$, $x_j = n-f$, $x_i = n-g$, $y_i = n+h$ for some positive integers $f$, $g$, $h$. Since $1 \le y_i - x_i, y_j - x_j \le 2c$, $f, g, h \le 2c$. We have $n (n - f) = (n - g)(n + h)$ which implies $(f + h - g) n = g h$. Since $g h > 0$, we must have $f + h - g > 0$. Therefore
\[
\sqrt{\frac{2N}{4 c^2}} \le \sqrt{\frac{2N}{\mu_i}} \le n \le (f + h - g) n = g h \le (2 c)^2
\]
which contradicts $N > 32 c^6$.

\section{Simultaneous Pell equations} \label{simulpell}
Summing up, if $N^2 = (N - d_1)(N + e_1) = (N - d_2)(N + e_2) = ... = (N - d_r)(N + e_r)$ with $1 \le d_1 < d_2 < ... < d_r \le c N^{1/2}$ and $1 \le e_1 < e_2 < ... < e_r \le c N^{1/2}$ and $N > 32 c^6$, then we have $2N = \mu_1 x_1 y_1 = \mu_2 x_2 y_2 = ... = \mu_r x_r y_r$ where $\mu_i$'s are distinct, $\mu_i \le 4c^2$ and $1 \le y_i - x_i \le 2c$. Let $y_i := x_i + c_i$ for some integer $1 \le c_i \le 2c$. Then
\[
8N = \mu_1 (2x_1) (2x_1 + 2c_1) = \mu_2 (2x_2) (2x_2 + 2c_2) = ... = \mu_r (2x_r) (2x_r + 2c_r).
\]
Suppose $r \ge 3$ for otherwise Theorem \ref{thm1} is true. Now $\mu_1 (2x_1 + c_1)^2 - \mu_1 c_1^2 = \mu_1 (2x_1) (2x_1 + 2c_1) = \mu_2 (2x_2) (2x_2 + 2c_2) = \mu_2 (2x_2 + c_2)^2 - \mu_2 c_2^2$. This leads to the Pell equation
\begin{equation} \label{pell1}
\mu_1 (2x_1 + c_1)^2 - \mu_2 (2x_2 + c_2)^2 = \mu_1 c_1^2 - \mu_2 c_2^2.
\end{equation}
Similarly,
\begin{equation} \label{pell2}
\mu_1 (2x_1 + c_1)^2 - \mu_3 (2x_3 + c_3)^2 = \mu_1 c_1^2 - \mu_3 c_3^2.
\end{equation}
\begin{lem} \label{lem1}
For $i \neq j$, $\mu_i c_i^2 \neq \mu_j c_j^2$.
\end{lem}
Proof: Suppose $\mu_i c_i^2 = \mu_j c_j^2$ for some $1 \le i < j \le r$. Then
\[
\frac{2x_i}{c_i} \Bigl( \frac{2x_i}{c_i} + 1 \Bigr) = \frac{\mu_i (2x_i) (2x_i + c_i)}{\mu_i c_i^2} = \frac{\mu_j (2x_j) (2x_j + c_j)}{\mu_j c_j^2} = \frac{2x_j}{c_j} \Bigl( \frac{2x_j}{c_j} + 1 \Bigr)
\]
which implies $\frac{x_i}{c_i} = \frac{x_j}{c_j} = \Lambda > 0$ say. Then
\[
\mu_i x_i^2 \Bigl(1 + \frac{1}{\Lambda} \Bigr) = \mu_i x_i (x_i + c_i) = \mu_j x_j (x_j + c_j) = \mu_j x_j^2 \Bigl(1 + \frac{1}{\Lambda} \Bigr).
\]
Hence $\mu_i x_i^2 = \mu_j x_j^2$. But recall $2(N - d_i) = \mu_i x_i^2$ and $2(N - d_j) = \mu_j x_j^2$. This implies $2(N - d_i) = 2(N - d_j)$ which contradicts $d_i < d_j$.

\bigskip

By Lemma \ref{lem1}, we have $\mu_1 (\mu_1 c_1^2 - \mu_3 c_3^2) \neq \mu_1 (\mu_1 c_1^2 - \mu_2 c_2^2)$. By a result of Turk [\ref{T}, Proposition 3], the solutions to both (\ref{pell1}) and (\ref{pell2}) satisfy
\[
2x_1 + c_1 < e^{C (4c^2)^2 (\log 4c^2)^3 (4c^2 \log 4c^2) \log(4c^2 \log 4c^2)} \le e^{C' c^6 (\log c)^5}
\]
for some large constants $C$ and $C'$. But $x_1 + c_1 = y_1 > \sqrt{\frac{2N}{\mu_1}} \ge \sqrt{\frac{2N}{4c^2}}$. This gives a contradiction if $N > e^{C'' c^6 (\log c)^5}$ with $C''$ sufficiently large. Therefore if $N^2 = (N - d_1)(N + e_1) = (N - d_2)(N + e_2) = ... = (N - d_r)(N + e_r)$ with $1 \le d_1 < d_2 < ... < d_r \le c N^{1/2}$ and $1 \le e_1 < e_2 < ... < e_r \le c N^{1/2}$ and $N > e^{C'' c^6 (\log c)^5}$, then $r \le 2$.

\section{A catch}
The above argument is almost correct except that when applying Turk's result to simultaneous Pell equations, one requires the coefficients $\mu_1$, $\mu_2$, $\mu_3$ in (\ref{pell1}) and (\ref{pell2}) to be squarefree. So we need to modify our argument. Suppose $\mu_1 = \tilde{\mu}_1 t_1^2$, $\mu_2 = \tilde{\mu}_2 t_2^2$, $\mu_3 = \tilde{\mu}_3 t_3^2$ where $t_i^2$ is the largest perfect square that divides $\mu_i$ and hence $\tilde{\mu}_i$ is squarefree. Then since $\mu_i \le 4c^2$, we have $\tilde{\mu}_i \le 4 c^2$ and $1 \le t_i \le 2 c$ for $i = 1,2,3$. The Pell equations (\ref{pell1}) and (\ref{pell2}) become
\begin{equation} \label{pell3}
\tilde{\mu}_1 [t_1 (2 x_1 + c_1)]^2 - \tilde{\mu}_2 [t_2 (2 x_2 + c_2)]^2 = \tilde{\mu}_1 t_1^2 c_1^2 - \tilde{\mu}_2 t_2^2 c_2^2,
\end{equation}
and
\begin{equation} \label{pell4}
\tilde{\mu}_1 [t_1 (2 x_1 + c_1)]^2 - \tilde{\mu}_3 [t_3 (2 x_3 + c_3)]^2 = \tilde{\mu}_1 t_1^2 c_1^2 - \tilde{\mu}_3 t_3^2 c_3^2.
\end{equation}
Turk's result requires $\tilde{\mu}_1 \neq \tilde{\mu}_2$ and $\tilde{\mu}_1 \neq \tilde{\mu}_3$. Suppose on the contrary $\tilde{\mu}_1 = \tilde{\mu}_2$ (the case $\tilde{\mu}_1 = \tilde{\mu}_3$ is similar). Tracing back, as $2N = \mu_1 x_1 y_1 = \mu_2 x_2 y_2$, we have
\[
2N = \tilde{\mu}_1 t_1^2 x_1 y_1 = \tilde{\mu}_2 t_2^2 x_2 y_2 \text{ which implies } (t_1 x_1) (t_1 y_1) = (t_2 x_2) (t_2 y_2).
\]
Note that $t_1 x_1 < t_1 y_1$ and $t_2 x_2 < t_2 y_2$. There are two possibilities.

\bigskip

Case 1: $t_1 x_1 = t_2 x_2$. Then we also have $t_1 y_1 = t_2 y_2$. Recall $y_i = x_i + c_i$. So $t_1 (x_1 + c_1) = t_2 (x_2 + c_2)$ which implies $t_1 c_1 = t_2 c_2$ and hence $\tilde{\mu}_1 t_1^2 c_1^2 = \tilde{\mu}_2 t_2^2 c_2^2$ or $\mu_1 c_1^2 = \mu_2 c_2^2$ which contradicts Lemma \ref{lem1}.

Case 2: $t_1 x_1 \neq t_2 x_2$. We are in the almost square situation in section \ref{almostsquare}. By the same argument except replacing $2c$ by $4c^2$ (since $t_i y_i - t_i x_i = t_i (y_i - x_i) \le 2c \cdot 2c$), we have a contradiction if $N > 512 c^{10}$.

\bigskip

Therefore, since in Theorem \ref{thm1} we have $n = N^2 > e^{C c^6 (\log c)^5}$, we do have $\tilde{\mu}_1 \neq \tilde{\mu}_2$ and $\tilde{\mu}_1 \neq \tilde{\mu}_3$. By Lemma \ref{lem1}, one can also check that $\tilde{\mu}_1 (\tilde{\mu}_1 t_1^2 c_1^2 - \tilde{\mu}_3 t_3^2 c_3^2) \neq \tilde{\mu}_1 (\tilde{\mu}_1 t_1^2 c_1^2 - \tilde{\mu}_2 t_2^2 c_2^2)$. Consequently, we can apply Turk's result to the Pell equations (\ref{pell3}) and (\ref{pell4}) and get the same result as in section \ref{simulpell} which finishes the proof of Theorem \ref{thm1}.

Department of Arts and Sciences \\
Victory University \\
255 N. Highland St., \\
Memphis, TN 38111 \\
U.S.A. \\
thchan@victory.edu

\end{document}